\documentclass[12pt]{amsart}
\numberwithin{equation}{section}

\begin{document}
\oddsidemargin=0cm
\evensidemargin=0cm


\title[Mappings between real submanifolds in
complex space]{ Mappings between real
submanifolds in complex space}
\author[Linda Preiss Rothschild]{ Linda Preiss
Rothschild}
\dedicatory{To Robert E. Greene on the occasion
of his 60th birthday}

\address{Department of Mathematics, 0112,
University of California at San Diego, La Jolla,
CA 92093-0112, USA}
\email{lrothschild@ucsd.edu}
\begin{thanks}{2000 {\em Mathematics Subject
Classification.} 32H02, 32V40,
32V35.}\end{thanks}
\begin{thanks}{The author is partially supported
by National Science Foundation grant
DMS-0100330.}\end{thanks}


\def\Label#1{\label{#1}}
\def\1#1{\overline{#1}}
\def\2#1{\widetilde{#1}}
\def\3#1{\mathcal{#1}}
\def\4#1{\widehat{#1}}

\def\bih{\sim}
\def\crd{\dim^{\rm{CR}}}
\def\crc{\codim^{\rm{CR}}}
\def\const{{\rm const}}

\def\s{s}
\def\k{\kappa}
\def\span{\text{\rm span}}
\def\tr{\text{\rm tr}}
\def\xo {{x_0}}
\def\Rk{\text{\rm Rk\,}}
\def\sg{\sigma}
\def\prt#1{{\partial \over\partial #1}}
\def\det{{\text{\rm det}}}
\def\po{p_0}
\def\fb {\bar f}
\def\gb {\bar g}
\def\Fb {\ov F}
\def\Gb {\ov G}
\def\Hb {\ov H}
\def\zb {\bar z}
\def\wb {\bar w}
\def \qb {\bar Q}
\def \t {\tau}
\def\z{\chi}
\def\w{\tau}
\def\Z{\zeta}
\def\g{\gamma}

\def \T {\theta}
\def \Th {\Theta}
\def \L {\Lambda}
\def\b {\beta}
\def\a {\alpha}
\def\o {\omega}
\def\l {\lambda}
\def \C{\mathbb C}
\def \R{\mathbb R}
\def \bR{\mathbb R}

\def\then{\Longrightarrow}
\def \im{\text{\rm Im }}
\def \re{\text{\rm Re }}
\def \Char{\text{\rm Char }}
\def \supp{\text{\rm supp }}
\def \codim{\text{\rm codim }}
\def \Ht{\text{\rm ht }}
\def \Dt{\text{\rm dt }}
\def \hO{\widehat{\mathcal O}}
\def \cl{\text{\rm cl }}
\def \bC{\mathbb C}
\def \C{\mathbb C}
\def \bL{\mathbb L}
\def \bZ{\mathbb Z}
\def \bN{\mathbb N}
\def \hol{\text{\rm hol}}
\def \aut{\text{\rm aut}}
\def \Aut{\text{\rm Aut}}
\def \J{\text{\rm Jac}}
\def\jet#1#2{J^{#1}_{#2}}
\def\gp#1{G^{#1}}
\def\gpo{\gp {2k_0}_0}
\def\emmp {\scrF(M,p;M',p')}
\def\rk{\text{\rm rk}}
\def\Orb{\text{\rm Orb\,}}
\def\Exp{\text{\rm Exp\,}}
\def\Span{\text{\rm span\,}}
\def\d{\partial}
\def\D{\3J}
\def\pr{{\rm pr}}
\def\Re{{\rm Re}\,}
\def\Im{{\rm Im}\,}

\newtheorem{Thm}{Theorem}[section]
\newtheorem{Def}[Thm]{Definition}
\newtheorem{Cor}[Thm]{Corollary}
\newtheorem{Pro}[Thm]{Proposition}
\newtheorem{Lem}[Thm]{Lemma}
\newtheorem{Conj}[Thm]{Conjecture}

\theoremstyle{remark}
\newtheorem{Exa}{Example}[section]
\newtheorem{Exas}[Exa]{Examples}
\newtheorem{Rem}[Exa]{Remark}

\maketitle

\begin{abstract}
In this paper I survey some recent results on finite
determination, convergence, and approximation of formal
mappings between real submanifolds in complex spaces.  A
number of conjectures are also given.
\end{abstract}

\section{Introduction}

In this survey we shall discuss a number of
questions and results on approximation,
convergence, and finite determination of
mappings between real submanifolds in complex
space.  We begin with some notation.
       Let $M$ and $M'$ be connected smooth real
submanifolds of $\Bbb C^N$ and $\Bbb C^{N'}$
respectively, with $p
\in M$ and
$p'
\in M'$. We write $H:(\C^N,p)\to (\C^{N'},p')$
for a germ  of a holomorphic mapping at
$p$  of $\C^N$ into $\C^{N'}$  with
$H(p) = p'$ .  We shall say that such a germ $H$
{\it sends
$M$ into
$M'$} if there is a neighborhood
$U$ of
$p$ in
$\C^N$ with
$H(M\cap U) \subset M'$.

       We shall consider here three questions
concerning mappings sending one real submanifold
into another.  The first question (and easiest
to state) is that of ``finite determination."
\medskip

\noindent {\bf Question 1.} Suppose
$H:(\C^N,p)\to (\C^{N'},p')$ is a germ of a
holomorphic mapping at $p$ sending $M$ into
$M'$.     Is there a positive integer $K$
(depending on $M$, $M'$, $p$,
$p'$ and possibly also on
$H$) such that if $\widetilde H$ is another such
germ of a holomorphic mapping  sending $M$ into
$M'$ and satisfying
\begin{equation}\Label {FDE}  \partial^\alpha
H(p) =
\partial^\alpha \widetilde H(p),  \
\
\
\forall  \alpha \in \Bbb N^N, \
       |\alpha| \le K,
\end{equation} then $H = \widetilde H$.
\medskip

It is easy to see that the answer to Question 1
is ``no" in general, even if $N=N'$, $M$ is real
analytic, and $H$ is a biholomorphism.  Indeed,
one may take $N=1$, $M = M'= \R$, and
$p=p'=0$.  In that case, any convergent power
series
$\sum_{n=1}^\infty a_nz^n$ with $a_n$ real for
all $n$ and
$a_1\not= 0$ defines a germ of a holomorphic
mapping at $0$ sending $M$ into $M'$.   However,
it follows from the ground breaking work of
Chern and Moser \cite{CM} in the mid '70's  that
the answer to Question 1 is ``yes" when $N=N'$,
$M$ and $M'$ are Levi nondegenerate
hypersurfaces, and $H$ is invertible i.e.\ the
Jacobian of $H$ is invertible at $p$.  In  fact,
in this case it is shown that one may take $K=2$.
        In this survey we shall give some recent
generalizations of this result and state some
open problems related to Question 1.
\medskip

In order to state two other  questions, we
introduce the notion of a formal mapping.  For
$p \in \C^N$ and
$p'\in \C^{N'}$,  a {\it formal mapping} $F:
(\C^N,p) \to (\C^{N'},p')$ is an $N'$-vector of
formal power series in
$N$ complex variables,
\begin{equation} F= (F_1,\ldots, F_{N'}), \ \ \
F_j(Z) =
\sum_{{\alpha\in \Bbb N}^N }a_{j\alpha}
(Z-p)^{\alpha},
\end{equation} with $(Z-p)^\alpha\colon =
(Z_1-p_1)^{\alpha_1}\ldots (Z_N-p_N)^{\alpha_N}$
and $F(p) = p'$.

        Since a formal mapping is not necessarily  a
holomorphic mapping, one must redefine what it
means to send one manifold into another. Suppose
  that  $(M,p)$ and $(M',p')$ are germs of
real analytic, real submanifolds in $\bC^N$ and
$\bC^{N'}$ respectively, given by the vanishing
of real-analytic (vector valued) real local
defining functions
$\rho(Z,
\bar Z)$ and
$\rho'(Z', \bar Z')$ near  $p\in M$ and $p'\in
M'$ respectively. We shall say that a formal
mapping $F$ as above sends
$M$ into
$M'$ if
$$\rho'\big(F(Z(x)),
\overline {F(Z(x))}\big) = 0$$ in the sense
of formal power series in $x$ for some (and
hence for any) real-analytic parametrization
$x\mapsto Z(x)$ of $M$ near $p=Z(0)$. (Here $x
\in \R^{\dim M}$.)

Our next question concerns convergence of formal
mappings.
\medskip

\noindent {\bf Question 2}.  Suppose
$F:(\C^N,p)\to (\C^{N'},p')$ is a formal mapping
taking
$M$ into
$M'$, where $M$ and $M'$ are real analytic
submanifolds of
$\Bbb C^N$ and
$\C^{N'}$ respectively.  Is $F$ necessarily
convergent?
\medskip

Again the case  $N=N' = 1$, $ M=M' = \R$, and
$p=p'=0$ is a simple counterexample.  Indeed, any
formal series $\sum^\infty_{n=1} a_n Z^n$  with
$a_n \in
\R$ for all $n$ represents a formal mapping
sending $M$ into $M'$.  Nevertheless, as in the
case of finite determination, it was proved in
\cite{CM} that if
$M$ and
$M'$ are Levi nondegenerate real analytic
hypersurfaces in
$\C^N$ and if
$F$ is  invertible (i.e. Jac$\, F(p) \not =0$),
then $F$ must be convergent.
\medskip

The third question we shall consider is more
subtle. Recall that one version of the Artin
approximation theorem \cite{Artin} implies that
if $F\colon (\C^j,0) \to (\C^N,0)$ is any formal
mapping and $h_\a\colon (\C^N,0) \to (\C^q,0)$ a
family of holomorphic mappings such that
$h_\a\circ F
= 0$ for all $\a$, then for any integer $k>0$
there is a holomorphic mapping $H^k\colon
(\C^j,0) \to (\C^N,0)$ such that the Taylor
series of $H^k$ agrees with that of $F$ up to
order $k$, and
$h_\a\circ H^k
= 0$ for all $\a$.  Question 3 asks whether
there is an analogous theorem for formal
mappings sending one submanifold into another.
\medskip

\noindent {\bf Question 3.}  Suppose
$F\colon(\C^N,p)\to (\C^{N'},p')$ is a formal
mapping sending
$M$ into
$M'$, where $M$ and $M'$ are real analytic, real
submanifolds of
$\Bbb C^N$ and
$\C^{N'}$ respectively. Can $F$ be approximated
by convergent mappings sending $M$ into $M'$?
More precisely,  given any $k > 0$ is there a
holomorphic mapping $H^k\colon (\C^N,p)\to
(\C^{N'},p')$ sending $M$ into $M'$ such that
\begin{equation}\Label {FDE2}  \partial^\alpha
H^{k}(p) =
\partial^\alpha F(p),  \ \ \
\forall  \alpha \in \Bbb N^N, \
       |\alpha| \le k ?
\end{equation}
\medskip

Returning to our example $N=N'
= 1$, $
M=M' = \R$, and
$p=p'=0$,   any formal mapping $F\colon(\C,0)\to
(\C,0)$ sending $\R$ into itself is given by
a power series with real coefficients.  Hence it
may be approximated by a convergent series to
any order by simply truncating the series.
Nevertheless, the answer to Question 3 is also
``no" in general.  Indeed, in their study of
real surfaces in
$\C^2$,  Moser and Webster \cite{MW} exhibited
pairs of real analytic surfaces
$M$ and $M'$ through $0$ for which there exists
a formal invertible mapping $F\colon(\C^2,0) \to
(\C^2,0)$ sending
$M$ into $M'$, but no germ of an invertible
holomorphic mapping
$H\colon(\C^2,0)
\to (\C^2,0)$ sending $M$ into $M'$.  In this
situation one says briefly that $M$ and $M'$ are
formally, but not biholomorphically, equivalent
at $0$.
\medskip

Note that Question 1, that of finite
determination, still makes sense if the
holomorphic mapping $H$ in that question is
assumed only to be formal. Before giving some
recent results, we will explore the relationship
between answers to Questions 1, 2, and 3.  It is
clear that a positive answer to Question 2  for
a fixed triplet ($(M,p)$, $(M',p')$, $F$), gives
a positive answer to Question 3 for that
triplet.  Now suppose that the answers to
Questions 1 and 3 are positive for  a fixed pair
of germs of real analytic, real submanifolds
$(M,p)$ and $(M',p')$ in $\C^N$ and for all
invertible formal mappings
      sending $p$ to $p'$ and $M$ to $M'$.
      We claim that this implies that  every such
mapping  is necessarily convergent.  Indeed, for
any invertible formal mapping $F$ sending $p$ to $p'$ and
$M$ to $M'$, by the positive answer to Question
1, there exists
$K > 0$ such that if $\widetilde F$ is another
such mapping
      and satisfies
\begin{equation}\Label {FDE3}  \partial^\alpha
F(p) =
\partial^\alpha \widetilde F(p),  \
\
\
\forall  \alpha \in \Bbb N^N, \
       |\alpha| \le K,
\end{equation} then $\widetilde F= F$.  By the
approximation of $F$ given by the positive
answer to Question 3 with $k=K$, there is a germ
of a holomorphic mapping $H\colon (\C^N,p)
\to (\C^{N'},p')$ sending $M$ into $M'$ and
satisfying (\ref{FDE3}) with $\2 F$ replaced by
$H$.  In that case it follows that
$H = F$ and hence $F$ is convergent.

In the rest of this paper I shall state some
recent
      results related to Questions 1, 2, and 3 and
also mention some open problems and conjectures.
In the last section I will discuss some aspects
of the techniques used in the proofs.

\section{Approximation}
       I shall begin with the following
approximation result, which is a partial answer
to Question 3.  If $F$ and $F'$ are formal mappings
at $p$ we shall
write $F-F' = O(|Z-p|^K)$ if the series
of the components of  $F$ and
$F'$ agree up to order
$K-1$ at $p$.

\begin{Thm}\Label{BRZ} {\rm(Baouendi, Rothschild,
\& Zaitsev \cite{BRZ})} Let $M\subset \bC^N$ be a
connected real analytic submanifold.  Then there
exists a closed proper real analytic subvariety
$V\subset M$ such that for every $p\in
M\setminus V$, every real analytic submanifold
$M'\subset
\bC^N$ with $\dim_{\R}M'=\dim_{\R}M$, every
$p'\in M'$,  every integer
$K>1$, and every invertible formal mapping
$F\colon (\C^N,p)\to (\C^N,p')$ sending $M$ into
$M'$, there exists a holomorphic mapping
$H\colon (\C^N,p)\to (\C^N,p')$ sending $M$ into
$M'$ with
$H(Z)-F(Z) = O(|Z-p|^K)$.
\end{Thm}

As a consequence of this approximation theorem,
one may conclude in particular that for all
points on a given submanifold outside a proper
real analytic subvariety, formal equivalence to
another submanifold is the same as local
biholomorphic equivalence.  In the absence of
counterexamples, one may conjecture the
following generalization of Theorem \ref{BRZ} in
which the assumptions that $M$ and $M'$ are
equidimensional submanifolds in the same complex
space and the invertibility of the formal map
$F$ are dropped.

\begin{Conj} \Label{BRZc}  Let $M\subset \bC^N$
be a connected real analytic submanifold.  Then
there exists a closed proper real analytic
subvariety $V\subset M$ such that for every
$p\in M\setminus V$, every positive integer
$N'$, every real analytic submanifold $M'\subset
\bC^{N'}$, every $p'\in M'$,  every integer
$K>1$, and every  formal mapping $F\colon
(\C^N,p)\to (\C^{N'},p')$ sending $M$ into $M'$,
there exists a holomorphic mapping  $H\colon
(\C^N,p)\to (\C^{N'},p')$ sending $M$ into $M'$
with
$H(Z)-F(Z) = O(|Z-p|^K)$.
\end{Conj}

A recent result in this direction was proved by
Meylan, Mir, and Zaitsev \cite{MMZ}. To state
their result we need to recall some definitions.
Let
$M\subset\C^N$ be a real submanifold of
codimension $d$ and
$p_0\in M$. Let
$\rho= (\rho_1, \ldots \rho_d)$ be a set of real
valued defining functions for $M$ near $p_0$ with
linearly independent differentials.  For
$p\in M$  near
$p_0$, let
\begin {equation}\Label{r1} r(p):=d-\dim\
\span_\C\left\{\rho_{j,Z}(p,\1p): 1\le j\le d
\right\}.
\end{equation} Here $\rho_{j,Z}=(\d\rho_j/\d
Z_1,\ldots, \d\rho_j/\d Z_N)\in\C^N$ denotes the
complex gradient of
$\rho_j$ with respect to
$Z=(Z_1,\ldots,Z_N)$. The point $p_0\in M$ is
called  {\em CR}  if the mapping
$p \mapsto r(p)$ is constant for $p$ in a
neighborhood of
$p_0$ in $M$. (It is easy  to see that if $M$ is
real analytic then the set of points at which
$M$ is not CR is a proper real analytic
subvariety of $M$.)
        The submanifold $M$ is called CR if it is CR
at all its points.  If $M$ is CR, the set
$\3 V \colon=T^{(0,1)}(\C^N)\cap \C TM$ of
antiholomorphic vectors tangent to $M$  form a
bundle, as does the set of holomorphic vectors
$\overline{ \3 V} \colon=T^{(1,0)}(\C)\cap
\C TM$. The bundle $\3 V$ is called the {\em CR
bundle} of $M$. The CR submanifold
$M$ is said to be of {\em finite type} at
$p$ (in the sense of Kohn \cite{K}, Bloom-Graham
\cite{BG}) if the span at $p$ of the Lie algebra
generated by the sections of the bundles
$\3 V$ and
$\overline{ \3 V}$ equals $\C T_pM$.

\begin{Thm}\Label{MMZ} {\rm(Meylan, Mir, \&
Zaitsev
\cite{MMZ})}  Let $M
\subset
\C^N$ be a real analytic CR submanifold, $p\in
M$ a point of finite type, and $K$ a positive
integer. Then for any real algebraic subset
$M'
\subset
\C^{N'}$, any $p' \in M'$, and any formal mapping
$F:(\C^{N},p)
\to (\C^{N'},p')$ sending $M$ into $M'$, there
exists a germ of a holomorphic mapping
$H:(\C^{N},p)
\to (\C^{N'},p')$ sending $M$ into $M'$ with
$H(Z)-F(Z) = O(|Z-p|^K)$.
\end{Thm}

Although the hypothesis that $M'$ be algebraic is
essential for the proof of Theorem \ref{MMZ} as
stated, one can conjecture the following.

\begin{Conj} The conclusion of Theorem 2.3 still
follows if $M'$ is assumed  to be a real
analytic set, rather than an algebraic
set.
\end{Conj}

Another direction in which Theorem
\ref{BRZ} might be generalized is to reduce the
exceptional subvariety
$V$ (at which mappings may not be approximable)
to the points at which $M$ is not CR.  Since the
only known examples of germs
$(M,p)$ and $(M',p')$ of real analytic
submanifolds of $\C^N$ which are formally
equivalent but not biholomorphically equivalent
occur in cases where $M$ and $M'$ are not CR at
$p$ and $p'$ respectively,  the following
conjecture seems reasonable.

\begin{Conj}\Label{BRZc2}  Let $M\subset \C^N$
be a real analytic CR submanifold.  Then for
any $p \in M$, any real analytic
real submanifold
$M'\subset \C^N$ with $\dim_\R M = \dim_\R M'$,
and
$p'\in M'$, any formal invertible mapping
$F:(\C^N,p)\to (\C^N,p')$ sending $M$ into $M'$
can be approximated by holomorphic ones, as in
the conclusion of Theorem \ref{BRZ}.
\end{Conj}

In fact, when $M$ is assumed to be of finite
type at $p$, Baouendi, Mir, and the author
\cite{BMR} have recently proved this conjecture.

\section {Convergence of formal mappings } \Label{converge}
In this section we shall discuss results and
conjectures concerning the convergence of formal
mappings which take one real analytic
submanifold into another.  That is, we shall
seek sufficient conditions for germs $(M,p)$ and
$(M',p')$ of real analytic submanifolds of
$\C^N$ and $\C^{N'}$ and classes $\3 F$ of
formal mappings sending $M$ into $M'$ which
guarantee that every element of $\3 F$ is
convergent, i.e.\ holomorphic.

To describe the results and conjectures, we
first recall that a real CR manifold $M$ is
called {\em generic} at
$p$ if the number
$r(p)$ defined in (\ref{r1}) is $0$. (Any real
hypersurface is a generic submanifold.) If
$M$ is CR but
$r:=r(p) > 0$, then after a local change of
holomorphic coordinates  in a neighborhood of
$p$, one can assume $p=0$ and write
$M = M_1 \times 0 \subset
\C^{N-{r}}\times \C^{r}$, where $M_1$ is a real
analytic generic submanifold of
$\C^{N-{r}}$. (See e.g.\ \cite{book} for
details.) In this case, as is shown in Example \ref{indep}
below, there are invertible formal mappings taking
$(M,p)$ into itself which are not convergent.   It is also
clear that
     near $p$, $M$ may be identified with the
generic submanifold $M_1$ of
$\C^{N-{r}}$.  For this reason we
shall often restrict our attention to generic
submanifolds.

\begin{Exa}\Label{indep} Let $M=M_1\times \{0\}\subset
\C^{N-s}\times\C^s$, where $s$ is a
positive integer and
$M_1$ is a real submanifold through $0$ in $\C^{N-s}$.
Let
$e_N$ be the last basis vector in $\C^N$.  Then the
function $F:Z \mapsto Z + g(Z_N) e_N$ maps $(M,0)$ to
itself and is invertible for any formal series $g(Z_N)$
satisfying $g(0) = 0$ and $g'(0) = 0$. If $g$ is not
convergent, neither is $F$.  Furthermore, since $g$ is
arbitrary,
$F$ is not determined by any finite number of its
derivatives at $0$
\end{Exa}

Before stating a theorem, I will introduce some
nondegeneracy conditions for a generic
submanifold.  Suppose first that $M$ is a real
hypersurface in $\C^N$ given near
$p \in M$ by $\rho(Z, \bar Z) = 0$, where $\rho$
is a real-valued function
       with nonvanishing differential.  As above,
let
$\rho_Z$ be the vector $(\d\rho/\d Z_1,\ldots,
\d\rho/\d Z_N)$ and let $L_1,\ldots, L_{N-1}$ be
a basis of
$(0,1)$ vector fields tangent to $M$ near $p$.
The hypersurface
$M$ is {\em Levi-nondegenerate at $p$} if and
only if
\begin{equation}\Label {LND}
\span_\C\{\rho_Z(p), L_1\rho_Z(p),\ldots,
L_{N-1}\rho_Z(p)\} =
\C^N.
\end{equation} We introduce a generalization of
Levi-nondegeneracy not only for hypersurfaces
but also for generic submanifolds of higher
codimension.   A generic submanifold $M$ of
codimension $d$ is called {\em finitely
nondegenerate} at
$p$ if there is an integer $k$ for which
\begin{equation}\Label {kND}
\span_\C\{\rho_{j,Z}(p), L^\alpha\rho_{j,Z}(p),
\alpha
\in
\Bbb N^{N-d}, 1 \le j \le d, \  |\alpha| \le k\}
=
\C^N,
\end{equation} where $\rho =
(\rho_1,\ldots,\rho_d)$ are defining functions
of $M$ near $p$,
$L^\alpha = L_1^{\alpha_1}\ldots
L_{d}^{\alpha_{N-d}}$,
$|\a| = \a_1+\ldots \a_{N-d}$.  Here, as in the
case of a hypersurface,
$L_1,\ldots, L_{N-d}$ is a basis of
$(0,1)$ vector fields tangent to $M$ near $p$. If
$k$ is the smallest integer such that (\ref{kND})
holds, then
$M$ is said to be {\em
$k$-nondegenerate} at $p$.

The case of Levi degenerate, but finitely
nondegenerate,  hypersurfaces was considered in
the joint work of
     Baouendi, Ebenfelt and the author
     \cite{Asian}. In that work it was proved that
any invertible formal mapping sending one
finitely nondegenerate hypersurface in $\C^N$,
$N > 1$, onto another is convergent.  For
hypersurfaces in
$\C^N, N
\ge 2$, the condition of finite nondegeneracy at
a point can be satisfied only if the
hypersurface is already of finite type at that
point, but for generic submanifolds of higher
codimension  the two conditions are independent.
(In $\C$ any real analytic hypersurface (i.e.\
curve) is finitely nondegenerate at all points
but not of finite type at any point, since there
are no nontrivial (0,1) vector fields tangent to
a curve.)   In
\cite{rational} the above authors extended their
results to prove that an invertible formal
mapping
$F:(\C^N,p)
\to (\C^N,p')$ sending a generic submanifold
$(M,p)$ onto $(M',p')$ is necessarily convergent
provided that $M$ is finitely nondegenerate and
of finite type at
$p$.

Another condition that plays an important role
for these questions is that of {\em holomorphic
nondegeneracy}, which is weaker than  finite
nondegeneracy.  A generic submanifold $M \subset
\C^N$ is {\em holomorphically nondegenerate} at a
point $p$ if there is no nontrivial holomorphic
vector field (with holomorphic coefficients)
tangent to $M$ in a neighborhood of $p$.  The
condition was first introduced for hypersurfaces
by N. Stanton \cite{S} for studying infinitesimal
automorphisms, and more generally later by the
author in joint work with Baouendi and Ebenfelt
\cite{acta}
   in their study of algebraicity of
mappings.  The following gives some basic
properties of holomorphic nondegeneracy. See
e.g. \cite{book} for proofs.

\begin{Pro}\Label{HND} Let
$M\subset
\C^N$ be a connected, real analytic, generic
submanifold, and
$p
\in M$.  The following are equivalent:
\begin{enumerate}
\item $M$ is holomorphically nondegenerate at
$p$;
\item $M$ is holomorphically nondegenerate at
all $q \in M$;
\item $M$ is finitely nondegenerate on a dense
subset of points.
\end{enumerate}
\end{Pro}

In light of Proposition \ref{HND}, for a
connected, real analytic, generic submanifold $M$, one may
assume that $M$ is either holomorphically nondegenerate
at no points, in which case we shall say that
$M$ is holomorphically degenerate, or at all
points, in which case we shall say that
$M$ is holomorphically nondegenerate. In the
first case, one may use the implicit function theorem to show
that there is a dense open set
$U
\subset M$ such for any
$p
\in U$, after a change of local holomorphic
coordinates in $\C^N$ near
$p$, one may assume $p=0$ and $M = M_1 \times
\C^s$, where
$M_1$ is a holomorphically nondegenerate generic
submanifold of
$\C^{N-s}$, with $s$ a positive integer (see
\cite{BRZ}). Thus if
$M$ is connected and holomorphically degenerate, near
all points $p \in U$, after a change of holomorphic
coordinates near $p$, one may assume that $p=0$ and $M$
is given by a defining function which is independent of
one of the coordinates, say
$Z_N$.  In this case we may construct
invertible formal mappings $F:(\C^N,p)\to
(\C^N,p)$ which send
$M$ into itself but are not convergent by using the
mappings given in Example
\ref{indep} above.

    One can
also show (see e.g.\ \cite{Asian}) that such
mappings as those of Example \ref{indep} exist at all
$p \in M$ in the holomorphically degenerate case.
Hence another necessary condition for positive answers to
either Questions 1 or 2 is that the target manifold
$M'$ be holomorphically nondegenerate.

Recall that a formal mapping $F$ from $(\C^N,p)$
to itself is called {\em finite} if the ideal generated
by its  components is of finite codimension in the
ring,  $\C[[Z-p]]$, of all formal power series in $\C^N$
at
$p$.  The following theorem, which was first
proved by Mir \cite{Mir-hyp} in the case of
hypersurfaces, is the most general to date for
finite formal mappings from $(\C^N,p)$ to $(\C^N,p')$.
\begin{Thm}\Label{BMR}{\rm(Baouendi, Mir, \&
Rothschild
\cite{BMR})} Let
$M$ and
$M'$ be connected real analytic generic
submanifolds of $\C^N$ of the the same dimension,
and
$p\in M$ a point of finite type.  Let $p'\in M'$
and suppose that
$M'$ is holomorphically nondegenerate.  Then any
formal finite mapping
$F:(\C^N,p)\to (\C^N,p')$ which sends
$M$ into $M'$ is convergent.
\end{Thm}

It is clear that even under the  hypotheses on
$M$, $M'$ and $p$ in Theorem \ref{BMR}, one must
still impose stringent conditions on a formal
mapping
$F$ in order to assure convergence, as the
following example shows.

\begin{Exa}\Label{hole}Let $M = M' \subset \C^3$
be the hypersurface through $p = p' = 0$ defined
by
\begin{equation} M\colon =\{Z=(Z_1,Z_2,Z_3)\in
\C^3:
\im Z_3 = |Z_1|^2 - |Z_2|^2\},
\end{equation} and let $F_1(Z): (\C^3,0) \to
(\C,0)$  be any  divergent formal mapping.
Then the formal mapping $F= (F_1,F_1,0)$ sends
$M$ into itself, but is not convergent.
\end{Exa}

In Example (\ref{hole}) the target hypersurface
$M'$ contains the holomorphic subvariety $V=\{Z:
Z_1=Z_2\}$ and the formal map $F$ sends $M$ into
$V$.  If neither of the hypersurfaces $M$ nor
$M'$ contains a nontrivial holomorphic subvariety
through
$0$, then both are of finite type at $0$ and
holomorphically nondegenerate.  In this case, it
is possible to drop all conditions on the formal
mapping
$F$, as is shown by the following theorem.

\begin{Thm}\Label{hyp}
{\rm (Baouendi, Ebenfelt, \& Rothschild
\cite{conv})} Let
$M$ and $M'$ be real analytic hypersurfaces
through
$0$ in $\C^N$,
$N \ge 2$, and suppose that neither $M$ nor $M'$
contains a nontrivial holomorphic variety through
$0$. Then any formal mapping $F:(\C^N,0) \to
(\C^N,0)$ sending $M$ into $M'$ is convergent.
\end{Thm}

Even for real analytic connected hypersurfaces in
$\C^2$,
     necessary and sufficient conditions to
guarantee the convergence of all invertible
formal mappings sending one into another are not
known. In this case, finite type at any point
implies holomorphic nondegeneracy. A reasonable
conjecture is the following.

\begin{Conj}\Label{hypc2} Let $M$ and $M'$ be
connected real analytic hypersurfaces through
$0$ in
$\C^2$, and suppose that $M$ is of finite type at
some point.  Then any invertible formal mapping
$F:(\C^2,0)\to(\C^2,0)$ sending $M$ into $M'$ is
convergent.
\end{Conj}

If $M$ is of finite type at $0$, then the
conclusion of Conjecture \ref{hypc2}, is true;
this has been proved in \cite {conv}.  Hence
the conjecture addresses the situation in which
$M$ is of finite type only at a dense set of
points and
$M$ is not of finite type at $0$.

       The situation for generic submanifolds in
complex spaces of different dimensions is much
more complicated. The following  result has been
proved for the case where the target manifold is
assumed to be real algebraic, i.e.\ contained in a real
algebraic subvariety of the same dimension.

\begin{Thm} \Label{MMZthm} {\rm (Meylan, Mir, \&
Zaitsev
\cite{MMZ})} Let $M \subset \C^N$ be a real
analytic generic submanifold and $p\in M$ a point
of finite type. If $M' \subset \C^{N'}$ is a real
algebraic set which contains no complex
subvariety, then for any $p' \in M'$,  any
formal mapping
$F:(\C^N,p)\to (\C^{N'},p')$ sending $M$ into
$M'$ is convergent.
\end{Thm}

It is reasonable to conjecture that the condition
of algebraicity in Theorem \ref{MMZthm} is
superfluous.  Hence the following.

\begin{Conj}\Label{MMZconj} If $M$ and $p$ are as
in Theorem \ref{MMZthm} and $M' \subset \C^{N'}$
is a real analytic set
containing no nontrivial holomorphic varieties,
then any formal mapping $F:(\C^N,p)\to
(\C^{N'},p')$ sending $M$ to $M'$ is convergent.
\end{Conj}

In fact Conjecture \ref{MMZconj} is still open even for
embeddings
into strictly pseudoconvex hypersurfaces (which
necessarily cannot contain nontrivial
holomorphic varieties). Recently, an advance for mappings between
strictly pseudoconvex hypersurfaces was proved by Mir.

\begin{Thm}\Label{Mir1}{\rm  (\cite{Mir})} Any
formal embedding sending a real-analytic strictly
pseudoconvex hypersurface
$M\subset \C^N$ into another real-analytic
strictly pseudoconvex hypersurface
$M'\subset \C^{N+1}$ is convergent.
\end{Thm}

In light of this result, the following special case of
Conjecture  \ref{MMZconj}
seems more accessible.

\begin{Conj}\Label{Mirconj} For any integers
$N,N'$, $2 \le N \le N'$,  any formal embedding
sending a real-analytic strictly pseudoconvex
hypersurface
$M\subset \C^N$ into another real-analytic
strictly pseudoconvex hypersurface
$M'\subset \C^{N'}$ is convergent.
\end{Conj}

Recall that for $N' = N$ the conclusion of
Conjecture \ref{Mirconj} holds by the results in
\cite{CM} mentioned above. Theorem
\ref{Mir1} above shows that the conjecture is true for $N'
= N+1$.

Theorem \ref{Mir1} is a consequence of a  more
general result for which the target is assumed
to be only Levi nondegenerate, rather than
strictly pseudoconvex. (In this case $M'$ may
contain a nontrivial holomorphic subvariety, or even be
foliated by complex subvarieties, as in Example \ref{hole}.)
To state the more general result we need another
definition.   A formal mapping
$F: (\C^N,p)\to (\C^{N'},p')$ sending
$M$ into $M'$ is called CR {\em transversal} if
$dF(p)(\C T_pM)\not
\subset \3 V' _{p'}\oplus \overline{\3 V'
_{p'}}$, where $\3 V'$ is the CR bundle of $M'$,
and $\3 V' _{p'}$ is the fiber of that bundle at
$p'$.  It is known that any formal embedding from
one strictly pseudoconvex hypersurface into
another is necessarily CR transversal (see Lamel
\cite{Lamel}, and Ebenfelt-Lamel \cite{EL}). Hence the
following result implies Theorem
\ref{Mir1}.

\begin{Thm}\Label{main3}{\rm (Mir \cite{Mir})}
Any formal CR transversal mapping sending a real
analytic Levi nondegenerate hypersurface
$M\subset \C^N$ into another such hypersurface
$M'\subset
\C^{N+1}$ is convergent.
\end{Thm}

When $M$ and $M'$ are not strictly pseudoconvex,
Theorem
\ref{main3} is sharp in the sense that
$N+1$ cannot be  replaced by any larger integer
$N'$, as shown by the following example (see
\cite{Lamel}).

\begin{Exa} Let $M \subset \C^N$ and $M'\subset
\C^{N+2}$ be given respectively by
\begin{equation} \Label{}
\im Z_N = \sum_{j=1}^{N-1} |Z_j|^2, \ \ \ \im
Z'_{N+2} = \sum_{j=1}^{N} |Z'_j|^2-|Z'_{N+1}|^2,
\end{equation} and consider the formal mapping
$F\colon(\C^N,0)\to (\C^{N+2},0)$ given by
$$F(Z)\colon
=(Z_1,\ldots,Z_{N-1},f(Z),f(Z),Z_N),$$ where
$f(Z)$ is a divergent series with $f(0)= 0$.  It
is easy to see that $F$ maps $M$ into
$M'$ and is CR transversal.  Here $M$ is
strictly pseudoconvex, but $M'$ is only Levi
nondegenerate and contains the complex curve given by
the equations
$Z'_{N+2}=Z'_1 =Z'_2= \ldots Z'_{N-1}=0$, $Z'_N =
Z'_{N+1}$.
\end{Exa}

For formal embeddings (not necessarily CR transversal)
for which the target manifold contains no nontrivial
complex curves, Mir has proved the following, which also
implies Theorem
\ref{Mir1}, since a strictly pseudoconvex hypersurface
cannot contain a complex curve.

\begin{Thm}\Label{main3}{ \cite{Mir}}
Any formal embedding sending a real
analytic Levi nondegenerate hypersurface
$M\subset \C^N$ into another such hypersurface
$M'\subset
\C^{N+1}$  is
convergent if $M'$ contains no nontrivial complex
curves.
\end{Thm}

\bigskip
\section{Finite determination of mappings:
results and conjectures}

In this section we return to variations of
Question 1.  As mentioned above, if $M$ and $M'$
are real analytic, one may ask the analogue of
       Question 1 when $H$ is assumed to  be a
formal mapping rather than a convergent one.
Since the publication of \cite{CM} in the mid '70's, a
number of results were obtained concerning
invertible mappings sending one finitely
nongenerate generic submanifold in $\C^N$ into
another. In
\cite{BER-cag} Baouendi, Ebenfelt, and the
author proved that if $M$ is a real analytic generic
submanifold of codimension
$d> 0$ and
$p\in M$ is a point of finite type at which $M$
is $k$-nondegenerate, then any germ at $p$ of an
invertible holomorphic mapping sending $M$ into
itself is determined by its derivatives at $p$
of order
$k(d+1)$.  This theorem has since been
generalized in several directions (see e.g.\
\cite{conv}), and a number of authors have also proved
stronger results in special cases.

For smooth CR submanifolds in $\C^N$ an
interesting class of mappings between two such
manifolds consists of the CR mappings, since
they are exactly the ones that preserve the CR
structures.  Recall that for CR submanifolds
$M\subset \C^N$, $M'\subset
\C^{N'}$, a differentiable mapping $h:M \to M'$
is called {\em CR} if
$h_*\3 V_p \subset \3 V'_{h(p)}$ for all $p \in
M$, where
$\3 V$ and $\3 V'$ denote the CR bundles of $M$
and
$M'$ respectively.  The restriction
of a holomorphic mapping is always CR, but when the
submanifolds are both real analytic, a smooth CR
mapping is the restriction of a holomorphic
mapping if and only if it is real analytic.

The following generalization to the CR case of
the theorem in \cite{BER-cag} cited above    was
proved first by Ebenfelt \cite{E-hyp} for the
case of hypersurfaces  and by Kim and Zaitsev
\cite{KZ} for higher codimension.

\begin{Thm} \Label{zk}{\rm (Ebenfelt
\cite{E-hyp}, Kim \& Zaitsev
\cite {KZ})} Let
$M, M'
\subset
\C^N$ be  smooth CR manifolds of codimension
$d$ with $M$
$k$-nondegenerate and of finite type at a point
$p$. If $h_1$ and $h_2$ are smooth invertible CR
mappings of $M$ into $M'$  such that
\begin{equation}\Label {FDEz}  \partial^\alpha
h_1(p) =
\partial^\alpha h_2(p),\ \ \forall \a,\ |\a| \le
k(d+1)
       \end{equation} then $h_1 = h_2$ in a
neighborhood of $p$.
\end{Thm}

In the case of formal or holomorphic mappings
there are some results in which the condition of
finite nondegeneracy in Theorem \ref{zk} can be
weakened. The notion of holomorphic nondegeneracy at a
point, as given in Section \ref{converge}, can be extended
to smooth manifolds by considering formal holomorphic vector
fields (see e.g.\cite{book}). For holomorphically
nondegenerate generic submanifolds the following was proved
in
\cite{BMR}.

\begin{Thm}\Label{fdhol} {\rm (Baouendi, Mir, \&
Rothschild \cite{BMR})} Let
$M\subset
\C^{N}$ be a smooth generic
submanifold through $p$, and assume that $M$ is of finite
type and holomorphically nondegenerate at
$p$. Then there exists a positive integer
$K$ such that for any smooth generic submanifold $M'
\subset \C^N$ with $\dim M = \dim M'$,  if
$H^j:(\C^N,p)\to (\C^{N},p')$, $j=1,2$, are
invertible formal mappings sending $M$ into $M'$
    and satisfying
\begin{equation}\Label {FDEz}  \partial^\alpha
H^1(p) =
\partial^\alpha H^2(p)\ \ \forall \a,\ |\a| \le
K,
       \end{equation}
       then $H^1=H^2$.
\end{Thm}

     For the case where the point
$p$ is assumed to be essentially finite (a
conditon weaker than finite nondegeneracy, but
stronger than holomorphic nondegeneracy, see
e.g.\
\cite {book}) this result was previously proved
in \cite{conv}.

For a hypersurface in $\C^2$ Ebenfelt, Lamel \&
Zaitsev \cite {ELZ}
have recently found optimal results for finite determination
of biholomorphisms.

\begin{Thm}\Label{sharp}
{\rm(\cite {ELZ})} Let
$(M,p)$ and
$(M',p')$ be germs of  real analytic
hypersurfaces in
$\C^2$ which are not Levi flat.  Then there is a
an integer
$K \ge 2$ such that if $H^1, H^2: (\C^2,p)\to
(\C^2,p')$ are germs of biholomorphisms sending
$M$ into
$M'$ with
\begin{equation}\Label {FDEz}  \partial^\alpha
H^1(p) =
\partial^\alpha H^2(p)\ \ \forall \a, \ |\a| \le
K,
       \end{equation}
       it follows that $H^1=H^2$. Moreover, if $M$
is of finite type at $p$, the conclusion holds
with $K = 2$.
\end{Thm}

Theorem \ref{sharp} is sharp.  Indeed, as
mentioned above, for a Levi flat hypersurface
(i.e.\ given by
$\im Z_2 = 0$), no finite $K$ exists for which
the conclusion of Theorem \ref{sharp} holds. In
addition, Kowalski \cite{kowalski} has given a
class of non-Levi flat (and hence holomorphically
nondegenerate) hypersurfaces
$(M,p)$ and
$(M',p')$ in $\C^2$ for which the best integer
$K$ as in Theorem \ref{sharp} can be arbitrarily
large.

       For hypersurfaces in $\C^N$, $N > 2$, finite
type no longer implies holomorphic nondegeneracy and is
not sufficient to guarantee  finite determination of
biholomorphisms, even between real analytic
hypersurfaces. Indeed, the mappings
$(Z_1,Z_2,Z_3) \mapsto (Z_1,Z_2,Z_3 +Z_3^K)$, $K
= 1, 2, \ldots$, are biholomorphisms at the
origin in
$\C^3$ sending the hypersurface
$M=\{(Z_1,Z_2,Z_3): \im Z_2 = |Z_1|^2\}$ into
itself; no finite $K$ will determine all such
mappings. Indeed, in this case the defining function of
$M$ is independent of $Z_3$ (so that $M$ is holomorphically
degenerate)  and the given mapping is a special case of
that of Example \ref{indep}.     Recently
it has been discovered that in $\C^N$, $N >2$, there are
real analytic hypersurfaces which are of finite type and
holomorphically nondegenerate at all points, but for which
self mappings fixing a point cannot be determined by
derivatives of order $2$.

  A survey including other results
on finite determination can be found in the recent article
by Zaitsev \cite{Z}.

\section{Segre mappings and the characerization
of finite type} In this last section we shall
define the so-called Segre mappings and give an
explicit example to show how these mappings can
be used to prove convergence of formal mappings
sending a generic submanifold into another. The closely
related Segre varieties (see below) have featured
prominently in the study of mappings such as in the work of
Webster
\cite{W} on algebraicity.  Segre mappings are also a
crucial tool in proving results concerning finite
determination and approximation.

For simplicity we shall take $M=M'\subset \C^N$
and $p=p' = 0$ in this section, where $M$ is a
real analytic generic submanifold of codimension
$d$ given near $0$ by the equation
$\rho(Z,\bar Z)=0$. Here
$\rho=(\rho_1,\ldots,\rho_d)$,   is a real
analytic local defining function of $M$ near $0$
with
$\partial_Z \rho_1(0) ,\ldots,\partial_Z
\rho_d(0)$ linearly independent.  We consider
$\bar Z$ as an independent variable so that
$\rho_j(Z,\zeta)$ is a convergent power series
in the
$2N$ indeterminates
$(Z,\zeta)$. We shall denote the ring of such
convergent power series with complex coefficients
by
$\bC\{ Z, \Z \}$.  The {\em generic rank}, $\Rk
F$, of  a formal mapping
$F\colon(\bC_x^k,0)\to(\bC_y^p,0)$   is defined
as the rank of the Jacobian  matrix $\d F /\d x$
regarded as a
$\Bbb K_x$-linear mapping
$\Bbb K_x^k\to\Bbb K_x^p$, where $\Bbb K_x$
denotes the field of fractions of $\bC[[ x]]$.
Hence
$\Rk F$ is the largest integer $s$ such that
there is an
$s\times s$ minor of the matrix $\d F/\d x$
which is not  0 as a formal power series in $x$.

Let $\gamma(\Z, t)$, where
$\Z=(\Z_1,\ldots\Z_N)$, $t=(t_1,\ldots,t_n)$, and
$n=N-d$, be a holomorphic mapping
$(\bC^N\times\bC^n,0)\to (\bC^N,0)$ such that
\begin{equation}\Label{gamma}
\rho(\gamma(\zeta,t),\zeta)= 0,\quad
\rk\, \frac{\d\gamma}{\d t}(0,0)=n.
\end{equation} The existence of such
$\gamma(\Z,t)$ is a consequence of the implicit
function theorem and the fact that
$\partial_Z\rho_1,\ldots,\partial_Z
\rho_d$ are linearly independent at $0$. We shall
call a holomorphic mapping $\gamma(\zeta,t)$
satisfying (\ref{gamma}) a {\it Segre variety
mapping} for the germ of $M$ at $0$.   The
mapping $t\mapsto
\gamma(\zeta,t)$, for $t$ near
$0\in\bC^n$,  parametrizes the {\em Segre
variety} of
$M$ at $\bar\Z$.

We define a sequence of holomorphic mappings
$v^j\colon (\bC^{nj},0)\to (\bC^N,0)$, called
{\it the iterated Segre mappings} of
$M$ at $0$ (relative to $\gamma$), inductively as
follows:
\begin{equation}
\begin{aligned}\Label{vj} &v^1(t^1):=
\gamma(0,t^1),\\ &v^{j+1}(t^1,\ldots,t^{j+1}):=
\gamma(\bar v^{j}(t^1,\ldots, t^j),t^{j+1}).
\end{aligned}
\end{equation} Here for $u(x) \in \C\{x\}$ we
denote by $\bar u(x)$ the series obtained by
replacing the coefficients in $u$ by their
complex conjugates. The relevance of the
iterated Segre mappings is given by the
following theorem.

\begin{Thm}\Label{main}{\rm (Baouendi, Ebenfelt,
\& Rothschild \cite{dyn})} Let
$M\subset
\bC^N$ be as above. Then
\begin{equation}\Label {idv}
\rho(v^{k+1}(t^1,\ldots,t^{k+1}),\bar
v^k(t^1,\ldots,t^{k}))= 0,\
\  k = 1,2,\ldots,
\end{equation}
\begin{equation}
\Rk\, v^1 \le \Rk\, v^2 \le \ldots \le \Rk v^j
\le
\Rk v^{j+1}
\le N, \ \ j = 1, 2, \ldots,
\end{equation} and
$M$ is of finite type at $0$ if and only if
$\Rk\, v^{d+1}=N$.  Moreover, for any $j$, $\Rk
v^j$ is independent of the choice of the Segre
variety mapping $\gamma$.
\end{Thm} We shall illustrate here how the Segre
mappings can be used to prove convergence and
finite determination in a very simple case.  We
let $M$ be the Lewy hypersurface in $\C^2$, i.e.\
\begin {equation} M = \{Z=(z,w) \in \C^2: \im w
= |z|^2\},
\end{equation} and let $\rho(Z,\Z) = w-\t -2i z
\chi$, where $\Z = (\chi,\t)\in \C^2$. Then
$\g=(\g_1,\g_2): (\C^3,0)\to (\C^2,0)$ is
uniquely defined by choosing $\g_1(\Z,t) = t$,
so that by (\ref{gamma}) and (\ref{vj}),
\begin{equation} \Label{}
\g_2(\Z,t) = \t + 2it\chi\ \ {\rm and}\ \
v^1(t^1) = (t^1,0),
\end{equation} Applying (\ref{vj}) for the next
iteration, we obtain
\begin{equation} \Label{} v^2(t^1,t^2) = \g(\bar
v^1(t^1),t^2) = (t^2,2it^1t^2), \ \  {\rm with}
\ \  \Rk v^2 = 2.
\end{equation} By Theorem \ref{main} we recover
the well-known fact, which is in this case
easier to prove directly from the definition,
that the Lewy hypersurface is of finite type at
$0$. For our purposes, however, it will be
necessary to go one further iteration to obtain
\begin{equation} \Label{} v^3(t^1,t^2,t^3) =
\g(\bar v^2(t^1,t^2),t^3) =
(t^3,-2it^1t^2+2it^2t^3).
\end{equation} Now suppose that $F: (\C^2,0) \to
(\C^2,0)$ is a formal invertible mapping sending
$M$ into itself.  Considering $\bar Z$ as an
independent variable $\Z$ it follows that
\begin{equation} \Label{} \rho(F(Z),\bar F(\Z))
= a(Z,\Z)\rho(Z,\Z),
\end{equation} for some $a(Z,\Z) \in
\C[[Z,\Z]]$.  Applying (\ref{idv}) of Theorem
\ref{main} with $k=2$ we obtain
\begin{equation} \Label{}
\rho(F(v^3(t^1,t^2,t^3)),\bar F(\bar
v^2(t^1,t^2)))= 0.
\end{equation}
    Writing
$F(z,w) = (f(z,w),g(z,w))$ and using our choice
of $\rho$, we hence have the identity
\begin{equation} \Label{meq}
g(t^3,2i(t^1t^2-t^2t^3)) - \bar g(t^2,2it^1t^2)
= 2if(t^3,2i(t^1t^2-t^2t^3))\bar
f(t^2,2it^1t^2)
\end{equation} Setting $t^1=t^2=0$ in
(\ref{meq}), and noting that
$f(0)=g(0)=0$, we obtain $g(z,0)= 0$. Since
$F$ is assumed to be invertible, it follows that
$f_z(0)\not=0$.  Hence after differentiating
(\ref{meq}) in $t^3$ we may solve  for
$\bar f(t^2,2it^1t^2)$ as a quotient of formal
series
\begin{equation} \Label{2}
\bar f(t^2,2it^1t^2) = \frac{
(g_z-2it^2g_w)(t^3,2i(t^1t^2-t^2t^3))}
{(f_z-2it^2f_w)(t^3,2i(t^1t^2-t^2t^3))}.
\end{equation} Now set $t^3=t^1$ in (\ref{2}) to
obtain
\begin{equation} \Label{2p1}
\bar f(v^2(t^1,t^2)) =\bar f(t^2,2it^1t^2) =
\frac{ g_z(t^1,0)-2it^2g_w(t^1,0)}
{f_z(t^1,0)-2it^2f_w(t^1,0)}=:R(t^1,t^2).
\end{equation}

To prove the convergence of the formal series
$F$ we shall make use of the the following
result (see e.g.\ \cite{EH} or
\cite{bierstone}).
\begin{Lem}\Label{conv} Let $S:(\C^p,0)\to \C$
be a formal series and $h:(\C^p,0)\to (\C^p,0)$
a germ of a holomorphic mapping with $\Rk h =
p$. Then if
$S\circ h$ is convergent, so is $S$.
\end{Lem}

       I claim that the convergence of $F$ will
follow from that of
$R(t^1,t^2)$. For this, assume first that
$R(t^1,t^2)$ is convergent. Since
$\Rk v^2 = 2$, it follows from Lemma
\ref{conv} that $\bar f(\chi,\t)$ (and hence
$f$) is also convergent. Setting $t^1=t^3$ in
(\ref{meq}), and making use of the fact that
$g(t^3,0)=0$, it follows from the convergence of
$f$ and again Lemma \ref{conv}  that
$\bar g(\chi,\t)$ is also convergent.  This
yields the claim that $F$ is convergent.

    It remains to show the convergence of
$R(t^1,t^2)$ given by (\ref{2p1}). By the
formula for
$R$, it will suffice to show that
\begin{equation} \Label{convt} t
\mapsto\partial^\a h(t,0) \ \ {\rm is\
convergent\ for}\ |\a| = 1,\  h = f\   {\rm
or}\  g.
\end{equation}

Putting $t^1=0$ in (\ref{2p1}) we obtain
\begin{equation} \Label{2p}
\bar f(t^2,0) = \frac{ g_z(0)-2it^2g_w(0)}
{f_z(0)-2it^2g_w(0)} =R(0,t^2),
\end{equation} where we note that $R(0,t^2)$ is
a rational (and hence convergent) function of
$t^2$ whose coefficients depend only on the first
derivatives of $F$ at
$0$. This shows that $\bar f(t^2,0)$ is
convergent.  Differentiating (\ref{2p}) we have
\begin{equation} \Label{2pd} \bar f_\chi(t^2,0) =
{d\over dt^2}R(0,t^2).
\end{equation} Since $g(z,0)= 0$, we have
$g_z(t,0) = 0$. This proves the convergence of
$f_z(t,0)$ and $g_z(t,0)$.

To complete the proof of the convergence of $F$,
it suffices to show that $t\mapsto \bar
f_\t(t,0)$ and
$t\mapsto \bar g_\t(t,0)$ are convergent.  For
the first, differentiate both sides of
(\ref{2p1}) with respect to
$t^1$ and set $t^1= 0$.  In fact this shows that
$\bar f_\t(t,0)$ is a rational function whose
coefficients depend only on the first two
derivatives of $f$ and $g$ at $0$. To prove
that  $t\mapsto g_w(t,0)$ is convergent, set
$t^3=t^1$ in (\ref{meq}) to obtain
\begin{equation} \Label{meq2}
       - \bar g(t^2,2it^1t^2) = 2if(t^1,0)\bar
f(t^2,2it^1t^2).
\end{equation} The desired convergence then
follows from that of
$\bar f(t^1,0)$ by differentiating (\ref{meq2})
with respect to $t^1$ and setting $t^1=0$.

\begin{Rem} In the above computation one can in
fact obtain an explicit formula for any local
biholomorphism mapping $M$ to itself, and it can
be readily seen that all such are given by
rational mappings.  Furthermore, it follows that
any such mapping is determined by its
derivatives of order two at $0$.
\end{Rem}

Although the above calculation is for a very
special case, an  argument along the same lines can be
used to show that if
$M\subset \C^N$ is a real analytic generic
submanifold which is finitely nondegenerate and
of finite type at
$p$, then any invertible formal mapping
$F:(\C^N,p) \to (\C^N,p)$ sending
$M$ into itself is convergent. (This is a special
case of Theorem \ref{BMR} stated above.)  In
fact, the proofs of most of the results stated
in this survey rely on Theorem
\ref{main} and also on a crucial use of Artin's
approximation theorem.


\end{document}